\documentclass[12pt]{amsart}
\usepackage{amssymb, amstext, amscd, amsmath}
\makeatletter
\def\@cite#1#2{{\m@th\upshape\bfseries%
[{#1\if@tempswa{\m@th\upshape\mdseries, #2}\fi}]}}
\makeatother
\theoremstyle{plain}
\newtheorem{thm}{Theorem}[section]
\newtheorem{cor}[thm]{Corollary}
\newtheorem{prop}[thm]{Proposition}
\newtheorem{lem}[thm]{Lemma}
\theoremstyle{definition}
\newtheorem{rem}[thm]{Remark}
\newtheorem{defn}[thm]{Definition}

\newcommand{\Prf}{\noindent\textbf{Proof.\ }}
\newcommand{\bx}{\hfill$\blacksquare$\medbreak}

\newcommand{\ca}{\mathrm{C}^*}

\newcommand{\ol}{\overline}

\newcommand{\wot}{\textsc{wot}}

\newcommand{\bbA}{{\mathbb{A}}}
\newcommand{\bbB}{{\mathbb{B}}}
\newcommand{\bbC}{{\mathbb{C}}}
\newcommand{\bbF}{{\mathbb{F}}}
\newcommand{\bbI}{{\mathbb{I}}}

\newcommand{\bbN}{{\mathbb{N}}}

\newcommand{\bbR}{{\mathbb{R}}}

\newcommand{\bbZ}{{\mathbb{Z}}}

  \newcommand{\A}{{\mathcal{A}}}
  \newcommand{\B}{{\mathcal{B}}}

  \newcommand{\E}{{\mathcal{E}}}

\renewcommand{\H}{{\mathcal{H}}}
  \newcommand{\I}{{\mathcal{I}}}
  \newcommand{\J}{{\mathcal{J}}}
  \newcommand{\K}{{\mathcal{K}}}

  \newcommand{\M}{{\mathcal{M}}}
  \newcommand{\N}{{\mathcal{N}}}

\renewcommand{\P}{{\mathcal{P}}}

\renewcommand{\S}{{\mathcal{S}}}
  \newcommand{\T}{{\mathcal{T}}}
  
  \newcommand{\V}{{\mathcal{V}}}

\newcommand{\fG}{{\mathfrak{G}}}

\newcommand{\fL}{{\mathfrak{L}}}
\newcommand{\fM}{{\mathfrak{M}}}

\newcommand{\qfor}{\quad\text{for}\quad}

\newcommand{\qforsome}{\quad\text{for some}\quad}

\newcommand{\Alg}{\operatorname{Alg}}
\newcommand{\alg}{\operatorname{alg}}

\newcommand{\rad}{\operatorname{rad}}

\newcommand{\spn}{\operatorname{span}}

\newcommand{\rep}{\operatorname{rep}}

\newcommand{\fgeeplus}{\bbF^+\!(G)}

\newcommand{\flgee}{\fL_G}

\begin{document}

\title[Isomorphisms of Algebras from Directed Graphs]{Isomorphisms of algebras associated
with directed graphs}
%
%
\author[E. Katsoulis]{Elias~Katsoulis}
\address{Department of Mathematics\\East Carolina University\\
Greenville, NC 27858\\USA}
\email{KatsoulisE@mail.ecu.edu}

\author[D.W. Kribs]{David~W.~Kribs}
\address{Department of Mathematics and Statistics\\University of Guelph\\
Guelph, Ontario\\CANADA N1G 2W1}
\email{dkribs@uoguelph.ca}
\begin{abstract}
Given countable directed graphs $G$ and  $G'$, we show that the
associated tensor algebras $\T_{+}(G)$ and $\T_{+}(G')$ are isomorphic as
Banach algebras if and only if the graphs $G$ are $G'$ are
isomorphic. For tensor algebras associated with graphs having no
sinks or no sources, the graph forms an invariant for algebraic
isomorphisms. We also show that given countable directed
graphs $G$, $G'$, the free semigroupoid algebras $\fL_G$ and $
\fL_{G'}$ are isomorphic as dual algebras if and only if the
graphs $G$ are $G'$ are isomorphic. In particular, spatially isomorphic free
semigroupoid algebras are unitarily isomorphic. For free
semigroupoid algebras associated with locally finite directed
graphs with no sinks, the graph forms an invariant for algebraic
isomorphisms as well.

\end{abstract}

\date{}
\maketitle

\section{ introduction and Preliminaries}   \label{S:intro}

Let $G$ be a countable directed graph with vertex set $\V(G)$, edge set
$\E(G)$ and range and source maps $r$ and $s$ respectively. The \textit{Toeplitz algebra
of $G$}, denoted as $\T(G)$, is the universal $\ca$-algebra generated by a set of partial
isometries $\{ S_e\}_{ e\in \E(G)}$ and projections $\{ P_x \}_{x\in \V(G)}$ satisfying
the relations
\[
(\dagger)  \left\{
\begin{array}{lll}
(1)  & P_x P_y = 0 & \mbox{for all $x,y \in \V(G)$, $ x \neq y$}  \\
(2) & S_{e}^{*}S_f = 0 & \mbox{for all $ e, f \in \E(G)$, $e \neq f $}  \\
(3) & S_{e}^{*}S_e = P_{s(e)} & \mbox{for all $e \in \E(G)$}      \\
(4)  & \sum_{r(e)=x}\, S_e S_{e}^{*} \leq P_{x} & \mbox{for all $x \in \V(G).$}
\end{array}
\right.
\]
The existence of such a universal object is implicit in \cite[Theorem 3.4]{Pim}
and \cite[Theorem 2.12]{MS2} and
was made explicit in \cite[Proposition 1.3 and Theorem 4.1]{FMR}.

\begin{defn}\label{tensor}
Given a countable directed graph $G$, the \textit{tensor algebra of $G$},
denoted as $\T_{+}(G)$, is the norm closed subalgebra of $\T(G)$
generated by the  partial
isometries $\{ S_e\}_{ e\in \E(G)}$ and projections $\{ P_x \}_{x\in \V(G)}.$
\end{defn}

The tensor algebras associated with graphs were introduced under the name \textit{quiver
algebras} by Muhly and Solel in \cite{MS2} and were further studied in \cite{MS1}.
They were defined rather
differently, in the setting of $\ca$-correspondences, but by \cite[Corollary 2.2]{FR}
our Definition~\ref{tensor} is equivalent to the original one.

The primary objective in this paper is to relate the structure of a tensor algebra to
its graph and show that the graph forms an invariant for bicontinuous
isomorphisms between tensor algebras of graphs (Theorem~\ref{main quiver}).
For tensor algebras associated with graphs having no
sinks or no sources, we show the graph forms an invariant for algebraic
isomorphisms as well (Corollary~\ref{quiveralgiso}).
For this purpose we prove an automatic continuity theorem, Theorem~\ref{automatic cntns},
of independent interest.

It is easy to see that
if two graphs $G$, $G'$ are isomorphic, then the tensor algebras $\T_{+}(G)$ and $\T_{+}(G')$
are isomorphic by the restriction of a
$\ca$-isomorphism between the associated Toeplitz algebras. Therefore the
focus will be on proving that if $\T_{+}(G)$ is isomorphic to $\T_{+}(G')$ as
a Banach algebra, then $G$ and $G'$ are isomorphic.
The proof of Theorem~\ref{main quiver} rests on analyzing certain representations
for a tensor algebra $\T_{+}(G)$. The one dimensional representations (characters)
of $\T_{+}(G)$ are parametrized by $\V(G)$ and points in the unit ball of a complex $n$-space for various $n$
(Proposition \ref{ball}). The two dimensional nest representations of $\T_{+}(G)$ allow us to identify the
edges in a way that reveals how they are "connected" and so reveals
the graph as an artifact of the representation theory of $\T_{+}(G)$ (Theorem~\ref{existence}). Independently,
similar ideas were used by Solel \cite{Sol} to prove an analogue of Theorem~\ref{main quiver}
for \textit{isometric} algebra isomorphisms. (The two papers overlap on the identification
of the character spaces, Proposition~\ref{ball} and Theorem~\ref{w*ball},
and the use of nest representations. Note however that the
representations in \cite{Sol} are assumed to be contractive while we are forced to
make no assumptions on ours, apart from continuity.) Our weaker hypothesis leads to
complications when "counting" edges between distinct vertices that need to be addressed
(Theorem~\ref{principal}).

Let $\lambda_{G,0}$ be the multiplication representation of $c_0 (\V(G))$ on $l^2 (\V (G))$, determined
by the counting
measure on $\V( G)$, and let $\lambda_{G}$ denote the representation of $\T(G)$ induced by $\lambda_{G,0}$,
in the sense of \cite{MS1} and \cite{FR}. It is easily seen that the Hilbert space $\H_G$ of
$\lambda_{G}$ is $l^{2}(\fgeeplus)$, where $\fgeeplus$ denotes the \textit{free semigroupoid} of the graph
$G$ (also called the \textit{path space} of $G$). This consists of all vertices $v \in \V(G)$ and all paths $w = e_k e_{k-1} \dots e_1$, where the
$e_i$ are edges satisfying $s(e_i ) = r( e_{i-1} )$, $i = 2, 3, \dots, k$, $k \in \bbN$. (Paths
of the form $w = e_k e_{k-1} \dots e_1$ are said to have length $k$, denoted as $|w|= k$, and vertices
are called paths of length $0$.) The maps $r$ and $s$ extend to $\fgeeplus$ in the obvious way,
two paths $w_1$ and $w_2$ are composable precisely when $s(w_2 ) = r( w_1 )$ and, in that case, the
composition $w_2 w_1$ is just the concatenation of $w_1$ and $w_2$. Let $\{ \xi_w \}_{w \in \fgeeplus}$
denote the usual orthonormal basis of $l^{2}(\fgeeplus)$, where $\xi_w$ is the characteristic function
of $\{ w \}$. Then, $\lambda_{G}(S_e )$, $e \in \E(G)$, equals the left creation
operator $L_e \in \B(l^{2}(\fgeeplus))$
defined by
\[
L_e \xi_w =
\left\{
\begin{array}{ll}
\xi_{ew} & \mbox{if $ s(e) = r(w)$} \\
0 & \mbox{otherwise.}
\end{array}
\right.
\]
(We shall write $P_v$ for $\lambda_{G}(P_v)$.) By \cite[Corollary 2.2]{FR},
the representation $\lambda_{G}$ is a faithful representation of $\T_{+}(G)$,
being the restriction of a faithful representation of the corresponding Toeplitz algebra.
For the rest of the paper, we will identify $\T_{+}(G)$ with its image under $\lambda_{G}$, i.e.,
we consider $\T_{+}(G)$ to be the norm closed
algebra generated by the operators $L_e$ and $P_v$ defined above.

\begin{defn}
The weak closure of $\lambda_{G}(\T_{+}(G)) \subseteq \B(l^{2}(\fgeeplus))$ is called the
\textit{free semigroupoid algebra} of $G$ and is denoted as $\flgee$.
\end{defn}

The second author and Power began a systematic study
for the free semigroupoid algebras in \cite{KrP, KP2}.
These algebras include as special
cases the space $H^{\infty}$ realized as the analytic Toeplitz
algebra \cite{Douglas_text,Hoffman_text} and the \textit{non
commutative analytic Toeplitz algebras} $\fL_n$, introduced by
Popescu \cite{Pop3, Pop2, Pop1} and  studied by Davidson and Pitts
\cite{DP1, DP2, DP3}, as well as Arias \cite{ArPop} and the
authors \cite{DKP, Kribs}. However, this is a broad enough class
to include various nest algebras, inflation algebras and infinite
matrix function algebras as special cases.

The second main result of the paper, Theorem~\ref{main freesem}, asserts that
given countable directed
graphs $G$ and $G'$, the free semigroupoid algebras $\fL_G$ and $
\fL_{G'}$ are isomorphic as dual algebras if and only if $G$ are $G'$ are isomorphic as graphs.
This improves one of the main results of
\cite{KrP}, which showed that the graph forms an invariant for unitary isomorphism between
free semigroupoid algebras, and provided the initial motivation for this paper.
Indepedently, Solel \cite{Sol} proved a similar result for \textit{isometric} $w*$-bicontinuous
isomorphisms. (Our result can be thought as an extension of Solel's result to Hardy algebras
\cite{MS3} associated with finite graphs.)
In particular, we show that spatially isomorphic free
semigroupoid algebras are unitarily isomorphic.
Our methods show, in fact, that with a slight additional hypothesis, the graph
forms an invariant for algebraic isomorphisms. While the arguments in Section~\ref{S:free sem}
follow the general line used in Section \ref{S:quiver}, the proofs are more involved.


\section{The classification of tensor algebras of graphs}    \label{S:quiver}

Given a countable directed graph $G$, the
character space of $\T_{+}(G)$ is denoted as $\fM_G$. It consists of all non-zero
multiplicative linear functionals on $\T_{+}(G)$ and is
equipped with the $w^*$-topology. Given
$x \in \V(G)$, the set of all $\rho \in \fM_G$ so that $\rho( P_x)=1$ is denoted as
$\fM_{G, x}$. Notice that the finite sums from
$\{P_x \}_{x \in \V(G)}$ form an approximate identity for $\T_{+}(G)$. Therefore, given $\rho \in
\fM_G$
there exists at least one $x \in \V(G) $ so that $\rho( P_x)=1$; the orthogonality
of $\{P_x \}_{x \in \V(G)}$ guarantees that such an $x$ is unique. Given any edge $e$ with
distinct source and range, the partial isometry $L _{e}$ is nilpotent
and so $\rho (L_e) = 0$, for any $\rho \in \fM_G$. Hence,
given $\rho \in \fM_{G,x}$ and a path $w$, we have $\rho (L_{w})=0$ whenever $w$ passes
through a vertex $y\neq x$.

\begin{prop} \label{ball}
Let $G$ be a countable directed graph. Let $x \in \V(G)$ and
assume that there exist exactly $n$ distinct loop edges whose
source is $x$. Then, $\fM_{G, x}$ equipped with the $w^*$-topology is
homeomorphic to the closed unit ball $\bbB_{n}$ inside $\bbC^n$.
\end{prop}

\Prf Let $e_1, e_2, \dots ,e_n$ be the distinct loop edges with source $x \in \V(G)$.
We define a map $\phi_x : \fM_{G, x} \longrightarrow \bbC^{n}$ by the formula
\[
\phi_x(\rho)= (\rho(L_{e_1}), \rho(L_{e_2}), \dots, \rho(L_{e_n})), \quad
\rho \in  \fM_{G, x}.
\]
Since any multiplicative form is completely contractive, and the $L_{e_i}$ are partial isometries
with mutually orthogonal ranges, it follows that the range of
$\phi_x$ is contained in $\bbB_n$. We will show that $\phi_x$ is the desired
homeomorphism.

Let $\lambda = (\lambda_1, \lambda_2, \dots \lambda_n)$ be an
n-tuple satisfying $\| \lambda \|_2^2 = \sum_{i=1}^{n}\, |
\lambda_i |^2 <1$; that is, $\lambda$ belongs to the interior of
$\bbB_{n}$. We define
\[
v_{\lambda , x}=(1 - \| \lambda \|_2^2) \sum_w \,
\overline{w(\lambda)} \xi_w,
\]
where $w$ in the above sum ranges over all monomials in $e_1, e_2,
\dots ,e_n$. (This construction originates from \cite{ArPop,DP2}
and was reiterated in \cite{KrP}). It is easily seen that
$v_{\lambda , x}$ is an eigenvalue for $\T_{+}(G)^{*}$ satisfying $P_x
v_{\lambda , x} = v_{\lambda , x}$ and $L_{e_i}^{*}v_{\lambda , x}
= \overline{\lambda_{i}}v_{\lambda , x}$. Therefore, the linear
form $\rho_{\lambda, x}$ defined as
\[
\rho_{\lambda, x}(A) = \langle Av_{\lambda , x} , v_{\lambda , x}\rangle ,
\quad A \in \T_{+}(G),
\]
belongs to $\fM_{G, x}$.

The above paragraph, combined with a compactness argument,
shows that $\phi_x (\fM_{G, x}) = \bbB_n$.
Now we show that $\phi_x$ is injective. If $\rho_1, \rho_2 \in \fM_{G,x}$ so that
$\phi_x ( \rho_1 ) = \phi_x ( \rho_2)$, then $\rho_1 (L_{e_i}) = \rho_2 (L_{e_i})$, $i=1,2,
\dots,n$.
As we have mentioned earlier, $\rho_1 (L_f) = \rho_2 (L_f) = 0$ for all other $f \in \E(G)$ and
so
$\rho_1$ and $\rho_2 $ agree on the generators of $\T_{+}(G)$. By continuity, $\rho_1 =\rho_2$.

We have established that $\phi_x$ is an continuous injective map between compact Hausdorff
spaces.
Hence $\phi_x^{-1}$ is also continuous.
\bx

\begin{rem}\label{fnlremark}
If $n =\infty$, then a easy
modification of the above proof shows that $\fM_{G, x}$ is
isomorphic to the closed unit ball of $l_2$, equipped with the
weak topology.
\end{rem}

\begin{cor}  \label{vertexcor}
If $G$ is a countable directed graph, then
$\fM_G$ is a locally compact Hausdorff space whose connected components coincide with
$\fM_{G, x}$, $x \in V(G)$.
\end{cor}

\Prf The space $\fM_G \cup \{0\}$ is a compact Hausdorff space and so $\fM_G$ is locally
compact. Fix an $x \in \V (G)$ and consider the map $\psi_x : \fM_{G}
\longrightarrow \bbC$ defined as
\[
\psi_x(\rho)=\rho (P_x), \quad \rho \in \fM_G.
\]
Clearly, $\psi_x$ is continuous and $\fM_{G , x} = \psi_x^{-1}(\{1\}) =  \psi_x^{-1}((0, 2))$.
Hence, $\fM_{G , x}$ is clopen. Proposition \ref{ball} shows that $\fM_{G , x}$ is
also connected and the conclusion follows.
\bx

The character space $\fM_G$ allows us to identify the vertices
of the graph $G$ with the connected components of $\fM_G$. In
order to decide whether there exists a directed edge between two
given vertices we use nest representation theory.

Let $\pi : \T_{+}(G) \longrightarrow \Alg \N$ be a two dimensional nest
representation of $\T_{+}(G)$, where $G$ is a countable directed graph.
Assume that the nest of projections $\N = \{0, N , I \}$ is acting
on a Hilbert space $\H$ and so there exists a basis $\{ h_1 , h_2
\}$ of $\H$ with $N = [h_2 ]$ and $N^\perp = [h_1]$. We may
associate with $\pi$ two multiplicative linear forms
$\rho_{\pi}^{( 1)}$ and $\rho_{\pi}^{( 2)}$, defined as
\begin{equation} \label{rhopi}
\rho_{\pi}^{( i)}(A) = \langle \pi(A) h_i , h_i \rangle, \quad A \in \T_{+}(G).
\end{equation}

\begin{lem} \label{edge}
Let $G$ be a countable directed graph, let $\pi : \T_{+}(G)
\longrightarrow \Alg \N$ be a two dimensional nest representation
and let $\rho_{\pi}^{( i)} \in \fM_{G, x_i}$, $i=1,2$, be as in
(\ref{rhopi}). If $x_1 \neq x_2$, then there exists an edge
$e \in \E(G)$ so that $e = x_2 e x_1$.
\end{lem}

\Prf By way of contradiction assume that there
are no edges in $G$ with source $x_1$ and range $x_2$.
There exist $a,b \in \bbC$ so
that
\[
\pi(P_{x_1}) =
\left(
\begin{matrix}
0 & a \\
0 & 1
\end{matrix}
\right)
\quad
\mbox{and}
\quad
\pi(P_{x_2}) =
\left(
\begin{matrix}
1 & b \\
0 & 0
\end{matrix}
\right).
\]
(However, $\pi(P_{x_2})\pi(P_{x_1})=0$ and this implies that $a = -b$.)

Let $S \in \T_{+}(G)$ so that $\pi(S) =
\left(
\begin{smallmatrix}
0 & 1\\
0 & 0
\end{smallmatrix}
\right)$. Without loss of generality we may assume that
$ S = P_{x_2}SP_{x_1} $ since
\[
\pi(P_{x_2}SP_{x_1})=
\left(
\begin{matrix}
1 & b \\
0 & 0
\end{matrix}
\right)
\left(
\begin{matrix}
0 & 1 \\
0 & 0
\end{matrix}
\right)
\left(
\begin{matrix}
0 & a \\
0 & 1
\end{matrix}
\right)
=
\left(
\begin{matrix}
0 & 1 \\
0 & 0
\end{matrix}
\right).
\]
Hence, $S$ can be approximated in norm by polynomials in $L_w$, with $w =
x_2 w x_1 \in \fgeeplus$. This implies the existence of one such $w=x_2 w x_1 \in \fgeeplus$
with $\pi(L_w) =
\left(
\begin{smallmatrix}
0 & a\\
0 & 0
\end{smallmatrix}
\right)$, for some non-zero $a \in \bbC$. However, by assumption in this case we have
\[
L_w = AL_e BL_f C,
\]
where $A,B, C \in \T_{+}(G)$ and $e, f \in \E(G)$ satisfy, $s(e) \neq r(e)$ and $s(f) \neq r(f)$.
Both $\pi(AL_e)$ and
$\pi( BL_f C)$ are then nilpotent of order two and so their product is zero.
So, $\pi(L_w) = 0$, which is a contradiction.
\bx

\begin{defn}
Given two vertices $x_1, x_2 \in \V(G)$  define
$\rep_{x_1,x_2}(\T_{+}(G))$ to be the collection of all two dimensional
nest representations $\pi:\T_{+}(G) \longrightarrow \Alg \N$ determined
as in Lemma~\ref{edge} by edges $e\in \E(G)$ with $e = x_2 e x_1$;
i.e., so that $\rho_{\pi}^{(i)} \in \fM_{G , x_i }$, $i=1,2$.
\end{defn}

\begin{thm}  \label{existence}
Let $G$ be a countable directed graph and let $x_1\neq x_2$ be two distinct vertices
in $\V(G)$. Then $\rep_{x_1,x_2}(\T_{+}(G))\neq \emptyset$
if and only if there exists an edge $e \in \E(G)$
with $s(e) = x_1$ and $r(e) = x_2$.
\end{thm}

\Prf One direction follows from Lemma~\ref{edge}. For the other
direction assume that there exists an edge $e = x_2ex_1$. Consider
the Hilbert space $\H_e \subseteq \H_G$ generated by the
orthogonal vectors $\xi_{x_1}, \xi_e$ and let $E_e$ be the
orthogonal projection on $\H_e$. Notice that for any edge $f \in
\E(G)$, $L_f^*\xi_{x_1} = 0$. In addition, either $L_f^*\xi_{e}=
\xi_{x_1}$  or $L_f^*\xi_{e}=0$, for all $f \in \E(G)$. In any
case, the space $\H_e$ is co-invariant by $\T_{+}(G)$ and so the
mapping $\pi_{e}(A)= E_e AE_e$, $A \in \T_{+}(G)$, defines a
representation for $\T_{+}(G)$ on $\H_e$.

Let $N$ be the subspace generated by $\xi_{e}$ and let $\N= \{ 0, N , I\}$.
The previous considerations show that
$N$ is invariant by $E_e \T_{+}(G) E_e$ and so $\pi_e$ maps into $\Alg \N$. In addition,
$\pi_e(P_{x_1})= N^{\perp}$, $\pi_e(P_{x_2})= N$ and
\[
\pi_e(L_e)=
\left(
\begin{matrix}
0 & 1 \\
0 & 0
\end{matrix}
\right).
\]
Hence, $\pi_e$ is surjective. Finally, it is clear that $\rho_{\pi_e}^{( i)} \in
\fM_{G , x_i }$, $i=1,2$, and so $\rep_{x_1,x_2}(\T_{+}(G))\neq \emptyset$.
\bx

In order to calculate the number of edges between two vertices $x, y$ of $G$, we need to
gain a better understanding of  $\rep_{x,y}(\T_{+}(G))$. We define
\[
\K_{x,y} =\bigcap\, \big\{ \ker\pi \mid \pi \in \rep_{x,y}(\T_{+}(G)) \big\}.
\]
Notice that each $\pi \in \rep_{x,y}(\T_{+}(G))$ induces a natural representation of $\T_{+}(G) / \K_{x,
y}$,
which we also denote as $\pi$, and is defined as
\[
\pi( A +  \K_{x, y} )= \pi(A), \quad A \in \T_{+}(G).
\]

Recall that the \textit{Jacobson radical} $\rad ( \A )$ of a complex algebra $\A$ is defined as the intersection of the kernels
of the irreducible representations of $\A$. It is known that for any complex algebra $\A$, the Jacobson radical $\rad (\A )$
coincides with the largest ideal of $\A$ consisting of quasinilpotent elements.
 
\begin{prop}   \label{radical}
Let $G$ be a countable directed graph and let $x, y \in \V(G)$ with $x \neq y$. Then
$A +  \K_{x, y} \in \rad (\T_{+}(G) / \K_{x,y} )$ if and only if $\pi ( A )^2 = 0$,
for all $\pi \in \rep_{x,y}(\T_{+}(G))$.
\end{prop}

\Prf Assume first that $A +  \K_{x, y} \in \rad (\T_{+}(G) / \K_{x,y} )$. Then given any
$\pi \in \rep_{x,y}(\T_{+}(G))$, $\pi (A)= \pi (A+  \K_{x, y})$ is quasinilpotent. Therefore,
$\pi (A)$ is a strictly upper triangular
$2\times 2$ matrix and so $\pi ( A )^2 = 0$.

Conversely, it is easily seen that the collection of all cosets $A +  \K_{x, y} \in \T_{+}(G) / \K_{x,y} $
satisfying $\pi ( A +  \K_{x, y})^2 = 0$, $\forall \pi \in \rep_{x,y}(\T_{+}(G))$, forms an ideal of
$\T_{+}(G) / \K_{x,y}$. Therefore, any such coset is contained in $\rad (\T_{+}(G) / \K_{x,y} )$ and the
conclusion follows.
\bx

\begin{lem} \label{reduction}
Let $G$ be a countable directed graph and let $x, y \in \V(G)$ with $x \neq y$. Then
\[
A +  \K_{x, y} = P_y A P_x +  \K_{x, y}
\]
for any $A +  \K_{x, y} \in \rad (\T_{+}(G) / \K_{x,y} )$.
\end{lem}

\Prf Let $\pi \in \rep_{x,y}(\T_{+}(G))$.
Proposition~\ref{radical}
implies the existence of a $c \in \bbC$ so that
$\pi(A)=
\left(
\begin{smallmatrix}
0 & c\\
0 & 0
\end{smallmatrix}
\right).$
Arguing as in the proof of
Theorem~\ref{edge} we obtain
\[
\pi(P_{y}A P_{x})=
\left(
\begin{matrix}
1 & b \\
0 & 0
\end{matrix}
\right)
\left(
\begin{matrix}
0 & c \\
0 & 0
\end{matrix}
\right)
\left(
\begin{matrix}
0 & a \\
0 & 1
\end{matrix}
\right)
=
\left(
\begin{matrix}
0 & c \\
0 & 0
\end{matrix}
\right).
\]
and so
\[
\pi(A-P_{y}A P_{x})=0.
\]
Since $\pi \in \rep_{x,y}(\T_{+}(G))$ was arbitrary, the conclusion follows.
\bx

\begin{lem} \label{Phi1}
Let $G$ be a countable directed graph and let $x, y \in \V(G)$ with $x \neq y$.
Assume that there exist
$n$ distinct edges $e_1 , e_2, \dots, e_n $ with source $x$ and range $y$. Let
$A = P_y A P_x  \in \K_{x,y}$ and let $A \sim \sum a_w L_w$ be the Fourier expansion of
$A$ \cite{KrP}. Then, $a_{e_i}= 0$, for all $i =1,2,\dots,n$.
\end{lem}

\Prf Let $\pi_{e_i} \in \rep_{x,y}(\T_{+}(G))$,
$i=1,2, \dots,n$ be as in the proof of Theorem~\ref{existence}.
As we have seen,
$
\pi_{e_i}(L_{e_i})=
\left(
\begin{smallmatrix}
0 & 1 \\
0 & 0
\end{smallmatrix}
\right),
$
while $\pi_{e_i}(L_f)=0$, for all other $f \in \E(G)$. Hence,
\[
\pi_{e_i}(A)= a_{e_i}
\left(
\begin{matrix}
0 & 1 \\
0 & 0
\end{matrix}
\right)
=0,
\]
and so $a_{e_i}=0$ for $i=1,2,\dots,n$.
\bx

Let $\I$ be a closed ideal of a Banach algebra $\B$ and let $\{ B_i\}_{i \in \bbI}$
be a subset of $\I$. The set $\{ B_i\}_{i \in \bbI}$ is said to be a \textit{generating set} for
$\I$ if the closed ideal generated by $\{ B_i\}_{i \in \bbI}$ equals $\I$. (In that case,
the $B_i$'s are said to be the \textit{generators} of $\I$.) The ideal $\I$ is said to be
\textit{$n$-generated} when $n$ is the smallest cardinality of a generating set. If
$n = \infty$, we understand $n$-generated to mean that there is no
finite generating set for $\I$.

\begin{thm}   \label{principal}
Let $G$ be a countable directed graph and let $x, y \in \V(G)$ with $x \neq y$.
Assume that there exist
$n$ distinct edges with source $x$ and range $y$. Then
$\rad (\T_{+}(G) / \K_{x,y} ) \subseteq \T_{+}(G) / \K_{x, y}$ is n-generated.
\end{thm}

\Prf Let $e_1, e_2, \dots, e_n$ be the edges with source $x$ and range $y$. First we show that
the collection,
\[
L_{e_{i}} + \K_{x, y}, \quad i=1, 2, \dots,n,
\]
is a generating set for $\rad (\T_{+}(G) / \K_{x,y} )$.

Indeed, let $A +  \K_{x, y} = P_y A P_x +  \K_{x, y} \in \rad (\T_{+}(G) / \K_{x,y} )$
by Lemma~\ref{reduction},
and let
\[
p(L) =  \sum_{w_k = yw_k x}\,a_{w_k} L_{w_k} \, + \, \K_{x,y}
\]
be a polynomial which is $\epsilon$-close to $A + \K_{x,y}$.
We may assume that each $L_{w_k}$ in the sum above factors as
$L_{w_k} = L_{u_k}L_{e_{i}}L_{v_k}$, for some $u_k, v_k \in \fgeeplus$
and $i\in \{1, 2, \dots,n\} $. (An $L_{w_k}$ not of this form is easily seen to belong to
$\K_{x,y}$ since $\pi (P_z) =0$ for $\pi \in \rep_{x,y}(\T_{+}(G))$ and $z\notin \{x,y\}$.)
Hence each $L_{w_k} + \K_{x, y}$ belongs to the ideal generated
by $L_{e_{i}} + \K_{x, y}$, $i=1, 2, \dots,n$, and so does $p(L)+ \K_{x, y}$.
But the polynomials
$p(L) + \K_{x, y}$
approximate $A +  \K_{x, y}$ and so
$L_{e_{i}} + \K_{x, y}$, $i=1, 2, \dots,n$, is a generating set for $\rad (\T_{+}(G) / \K_{x,y} )$.

By way of contradiction assume that there exists a generating set
\[
A_j + \K_{x, y}, \quad j=1, 2, \dots,m,
\]
for $\rad (\T_{+}(G) / \K_{x,y} )$ with $m<n$. By Lemma~\ref{reduction}, we may assume that
$A_j = P_y A_j P_x $. Therefore, there exist scalars $a_{j, i}$ so that
\[
\Phi_1(A_j)= a_{j,1}L_{e_{1}} +a_{j,2}L_{e_{2}}+\dots+a_{j,n}L_{e_{n}} ,
\]
where $\Phi_1$ is the contractive idempotent on $\T_{+}(G)$ defined by the formula
\[
\Phi_1\Big(\sum_w \, a_wL_w\Big)=\sum_{|w|=1} \, a_wL_w.
\]
Let $\M$ be the subspace of $\T_{+}(G)$ generated by $L_{e_{i}}$, $i=1,2,\dots,n$, and let
$\M_0 \subset \M$ be the linear span of
\[
a_{j,1}L_{e_{1}} +a_{j,2}L_{e_{2}}+\dots+a_{j,n}L_{e_{n}}, \quad j=1,2,\dots,m.
\]
Let $X \in \M \backslash \M_0$. Clearly, $X+ \K_{x,y}$ belongs to $\rad (\T_{+}(G) / \K_{x,y} )$.
Since the collection $A_j + \K_{x, y}$, $j=1, 2, \dots,m$, is generating
for $\rad (\T_{+}(G) / \K_{x,y} )$, the operator $X$ can be approximated by finite sums of the form
\begin{equation} \label{shorter}
\sum_{k}\, B_k \hat{A}_k C_k + D ,
\end{equation}
where, $B_k , C_k \in \T_{+}(G)$, and $D = P_y D P_x  \in \K_{x, y}$ and $\hat{A}_k \in
\{A_1, A_2, \dots A_m \}$. Lemma \ref{Phi1}
shows now that
$\Phi_1(D)= 0$. It is also easy to see that
\[
\Phi_1 (B A_j C) = \mu_j \big( a_{j,1}L_{e_{1}} +a_{j,2}L_{e_{2}}
+\dots+a_{j,n}L_{e_{n}}\big), \quad \mu_j \in \bbC,
\]
for all $j= 1,2,\dots,m$ and $B,C \in \T_{+}(G)$. Hence,
\begin{equation} \label{sum}
\Phi_1 \left( \sum_{k}\, B_k \hat{A}_k C_k + D \right) \in \M_0.
\end{equation}
A simple approximation argument with the contractive map $\Phi_1$ shows now
that $\Phi_1(X) = X$ belongs to $ \M_0$, which is a contradiction.
\bx

We are in position now to prove the classification theorem for the tensor
algebras of graphs. Given a directed graph $G$, we associate a graph $\fG$, which is constructed
using the algebraic structure of $\T_{+}(G)$ as follows. The vertices of $\fG$ are the connected
components of the character space $\fM_G$ of $\T_{+}(G)$. By Corollary \ref{vertexcor}
there is a natural one-to-one correspondence $g:\V(\fG) \longrightarrow\V(G)$
so that $\rho(P_{g(x)})=1$, for all $\rho \in x$. To each vertex $x \in\V(\fG) $, we attach
as many loop edges as the homeomorphic class of $x$. By Proposition~\ref{ball}, the
number of loop edges starting at $x$ and $g(x)$
are equal. Given two distinct vertices $x_1, x_2 \in\V(\fG)$, we create as many directed edges
from
$x_1$ to $x_2$ as the least number of generators for $\rad (\T_{+}(G) / \K_{g(x_1 ) , g(x_2 )} )
\subseteq \T_{+}(G) / \K_{g(x_1), g(x_2)}$, provided
that $\rep_{g(x_1),g(x_2)}(G) \neq \emptyset$. Theorem~\ref{principal}
shows that the number of directed edges from $x_1$ to $x_2$ and
$g(x_1)$ to $g(x_2)$ coincide. Hence, $G$ and $\fG$ are isomorphic as graphs.

\begin{thm} \label{main quiver}
Let $G$, $G'$ be countable directed graphs. The algebras $\T_{+}(G)$,  $ \T_{+}(G')$ are
isomorphic as Banach algebras if and only if the graphs $G$ are $G'$
are isomorphic.
\end{thm}

\Prf Assume that the algebras $\T_{+}(G)$,  $ \T_{+}(G')$ are
isomorphic as Banach algebras. It suffices to show that the graphs $\fG$ and $\fG^{\prime}$ are
isomorphic. Let
\[
\tau: \T_{+}(G) \longrightarrow \T_{+}(G')
\]
be a bicontinuous isomorphism between $\T_{+}(G)$ and $ \T_{+}(G')$. Then, $\tau$ induces
a homeomorphism between the
character spaces of $\T_{+}(G)$ and $ \T_{+}(G')$, which we still denote as $\tau$, and is defined
by the formula $\tau(\rho) = \rho \circ \tau^{-1}$, $\rho \in \fM_G$. Since homeomorphisms
preserve
connected components, $\tau$ establishes a one-to-one correspondence between the vertices
of $\fG$ and $\fG^{\prime}$, which we, once again, denote as $\tau$. It is clear that the number
of loop edges attached between $x$ and $\tau(x)$, $x \in \V(\fG)$, coincide.

Let  $x_1,x_2$ be two distinct connected components of $\fM_G$.
{\medbreak}
\textbf{Claim:} $\pi \in \rep_{x_1 , x_2}(G)$ if and only if $\pi\circ \tau^{-1} \in
\rep_{\tau(x_1) ,
\tau(x_2)}(G ')$.
{\medbreak}
Assume that $\pi \in \rep_{x_1 , x_2}(G)$.
Let $\N$ be a nest such that $\pi$ maps $\T_{+}(G)$ onto $\Alg \N$ and let $h_1, h_2$ be
as in (\ref{rhopi}) so that
\[
\rho_{\pi}^{( i)}(A) = \langle \pi(A) h_i , h_i \rangle, \quad A \in \T_{+}(G).
\]
Assume that, $\rho_{\pi}^{( i)} \in x_i$, $i=1,2$. Now $\pi \circ \tau^{-1}$
is a representation of $\T_{+}(G')$ and by definition,
\[
\rho_{\pi \circ \tau^{-1}}^{( i)}(A') = \langle \pi \circ \tau^{-1}(A') h_i , h_i \rangle, \quad A' \in
\T_{+}(G').
\]
Therefore, $\rho_{\pi \circ \tau^{-1}}^{( i)}(\tau(A))=  \rho_{\pi}^{( i)}(A)$ for $A\in \T_{+}(G)$
and so
\[
\rho_{\pi \circ \tau^{-1}}^{( i)} =  \rho_{\pi}^{( i)} \circ \tau^{-1} = \tau(\rho_{\pi}^{( i)}) \in
\tau(x_i)
\]
This proves one direction of the claim. By reversing the above argument, we obtain the other
direction.
{\medbreak}
The claim above implies now that $\tau (\K_{x_1 , x_2})= \K_{\tau(x_1) ,\tau( x_2)}$
and so $\tau$ induces an
isomorphism between $\T_{+}(G) / \K_{x, y}$ and $\T_{+}(G') / \K_{\tau(x_1) ,\tau( x_2)}$.
Furthermore,
\[
\tau\left( \rad (\T_{+}(G) / \K_{x_1 , x_2}) \right) = \rad (\T_{+}(G') / \K_{\tau(x_1) ,\tau( x_2)})
\]
and so $\rad (\T_{+}(G) / \K_{x_1 , x_2}) $ is $n$-generated
whenever $\rad (\T_{+}(G') / \K_{\tau(x_1) ,\tau( x_2)})$ is.
Hence $\tau$ preserves the number of directed edges between distinct
vertices of the graphs $\fG$ and $\fG^{\prime}$. This concludes the proof.
\bx

Given an (algebraic) isomorphism or an epimorphism between two Banach algebras, it is always
desirable to know
whether or not it is continuous. This is the problem of \textit{automatic continuity} and it has
attracted much attention since the early stages of Banach algebra theory. An important tool for
the study
of
this problem is Rickart's notion of a separating space. Let $\phi : \A \longrightarrow \B$ be an
epimorphism between
Banach algebras and let $\S (\phi )$ be  the \textit{separating space} of $\phi$. This is the
two-sided closed ideal of  $\B$ defined as
\[
\S (\phi ) = \big\{ b \in \B \mid \exists \{a_n\}_n \subseteq \A
\mbox{ such that } a_n \rightarrow 0 \mbox{ and } \phi(a_n )\rightarrow b\big\}.
\]
One can easily see that the graph of $\phi$ is closed if and only if $\S (\phi )= \{0\}$.
Thus by the closed graph theorem, $\phi$ is continuous
if and only if $\S (\phi ) = \{0\}$.

The following is an adaption of \cite[Lemma 2.1]{Sinclair} and was used in
\cite{DHK} for the study of isomorphisms between limit algebras.

\begin{lem}[Sinclair]     \label{Sinclair}
 Let  $\phi : \A \longrightarrow \B$ be an epimorphism between
Banach algebras and let $ \{B_n\}_{n\in \bbN}$ be any sequence in
$\B$. Then there exists $n_0 \in \bbN$ so that for all $n \geq n_0$,
\[
\overline{ B_1 B_2 \dots B_n \S (\phi ) } = \overline{ B_1 B_2 \dots B_{n+1}  \S (\phi ) }
\]
and
\[
\overline{ \S (\phi ) B_{n} B_{n-1} \dots B_1 } = \overline{\S (\phi ) B_{n+1} B_n \dots B_1 }.
\]
\end{lem}

The above lemma will allow us to classify a large class of tensor algebras up to algebraic
isomorphism.

\begin{defn}
A vertex $x$ of a countable directed graph is said to be a \textit{sink} (resp. \textit{source})
if it emits (resp. receives) no edges.
\end{defn}

\begin{thm}  \label{automatic cntns}
Let  $G$ be a countable directed graph which has either no sinks or no sources.
If $\A$ is any Banach algebra and $\phi : \A \longrightarrow \T_{+}(G)$ is an epimorphism of algebras
then $\phi$ is automatically continuous.
\end{thm}

\Prf  Consider the extension
$\hat{\phi}:\A + \bbC I\longrightarrow \T_{+}(G) + \bbC I$ of $\phi$.
It suffices to show that $\hat{\phi}$ is continuous.

Assume that $G$ has no sinks.
By way of contradiction assume that $\hat{\phi}$ is not continuous and so
$\S (\hat{\phi} )\neq \{0\}$. If there are infinitely many
distinct
$x_n \in \V(G)$, $n\in \bbN$, so that $P_{x_n} \S(\hat{\phi}) \neq 0$, then the sequence
\[
B_n = I-\sum_{i=1}^{n} \, P_{x_i}
\]
contradicts Lemma~\ref{Sinclair}. Hence there exists a vertex $x \in \V(G)$ and a path
$w=ywx$,
$y \in \V(G)$, so that $P_x  \S(\hat{\phi}) \neq 0$ and the vertex $y$ supports at least one loop,
say
$u \in \fgeeplus$. (Indeed, otherwise we could start from $x$, move forward on an infinite path with no
loops, by multiplying on the left with the corresponding creation operators,
and therefore produce infinitely many distinct $P_x$ with  $P_x \S(\phi) \neq 0$.)
Hence $P_y \S(\hat{\phi}) \neq 0$ and we have
\[
\overline{(L_u)^{n} \S (\hat{\phi} )} = (L_u)^{n} \S (\hat{\phi} )
\supset
(L_u)^{n+1} \S (\hat{\phi} ) = \overline{(L_u)^{n+1} \S (\hat{\phi} )}.
\]
This shows that the sequence $B_n = (L_u)^{n}$, $n \geq 1$,
contradicts Lemma~\ref{Sinclair} and so $\phi$ is necessarily
continuous.

An argument similar to the one above shows that $\hat{\phi}$ is also continuous when $G$ has
no sources.
\bx

\begin{cor}\label{quiveralgiso}
Let $G$, $G'$ be countable directed graphs which which have
no sinks or no sources. Then $\T_{+}(G)$ and
$ \T_{+}(G')$ are
algebraically isomorphic if and only if the graphs $G$ and $G'$ are isomorphic.
\end{cor}

In order to show that the graph of any tensor  algebra is an invariant for algebraic isomorphisms, one
has to extend Theorem~\ref{automatic cntns} to
arbitrary countable graphs. At the present, we do not know how to do this.
Nevertheless, there are cases where the
automatic continouity of the algebraic isomorphism is not
needed, as the following result shows.

\begin{cor}
Let $G$, $G'$ be countable directed graphs with no loop edges. Then $\T_{+}(G)$ and
$ \T_{+}(G')$ are
algebraically isomorphic if and only if the graphs $G$ and $G'$ are isomorphic.
\end{cor}

\Prf In this case we assume that $\rep_{x, y}(G)$ consists of not
necessarily continuous representations and we proceed as earlier. The
difference here is that $\rad (A_G / \K_{x,y})$ is an $n$-dimensional vector space,
provided that there exist exactly $n$-directed edges from $x$ to $y$. This property is preserved by
arbitrary isomorphisms and the conclusion follows.
\bx



\section{An Application to Free Semigroupoid Algebra Theory}  \label{S:free sem}

The free semigroupoid
algebras have been classified by Kribs and Power \cite{KrP}
up to unitary equivalence. Once again, the graph forms a complete
invariant. In light of Theorem~\ref{main quiver}, it is natural to
ask if the graph forms a complete invariant for bicontinuous or
even algebraic isomorphisms. In Theorem~\ref{main freesem} we
show that two free semigroupoid algebras are isomorphic as dual
algebras if and only if the associated graphs are isomorphic. In
particular, our result classifies the free semigroupoid algebras
up to similarity and shows that two such algebras are similar if
and only if they are unitarily equivalent. For free semigroupoid
algebras associated with locally finite directed graphs with no
sinks, we show that the graph forms an invariant for algebraic
isomorphisms as well.

The proof of Theorem~\ref{main freesem} follows the same line of
reasoning as that  of Theorem~\ref{main quiver}. All the results
of Section \ref{S:quiver}, as well as their proofs, adopt easily
to the $w^*$ context. The only exceptions are Proposition~\ref{ball}
and Theorem~\ref{principal}.

Let $\fM_{G}^{w^*}$ denote the set of all $w^*$-continuous
multiplicative linear functionals on $\fL_G$.
If $x \in \V(G)$ then $\fM_{G, x}^{w^*}$ denotes the
collection of all functionals $\rho \in \fM_{G}^{w^*}$ so that
$\rho (P_x) = 1$. Clearly, the (disjoint) union of all $\fM_{G,
x}^{w^*}$, $x \in \V(G)$, equals $\fM_{G}^{w^*}$.

\begin{defn}
Let $G$ be a countable directed graph and let $x \in \V(G)$.
We say that a directed loop $u = xux$ is \textit{primitive}
 if it does not factor as
$u = wv$, where $w, v \in \fgeeplus$, $r(w) =s(w)= r(v) =s(v) =x$,  and $1 \leq |w| < |u|$.
Clearly any loop edge supported at $x$ is a primitive loop but there may be many more. The
collection of
all primitive loops $w$ with $r(x) = s(x) =x$ is denoted as $\P (\fgeeplus, x)$.
\end{defn}

\begin{lem}  \label{secondunitequiv}
Let $G$ be a countable directed graph, $x \in \V (G)$ and let
$u_1, u_2 , \dots ,u_n$ , with $n= \infty $ possibly, be the
primitive loops with source $x$. Then there exists a
$w^*$-bicontinuous isomorphism from $P_x \fL_G P_x$ onto $\fL_n$,
which maps $L_{u_i}$ to $L_i$ for $i=1,2,\dots,n$.
\end{lem}

\Prf It suffices to show
that $P_x \fL_G P_x$ is unitarily equivalent to an ampliation of
$\fL_n$, such that generators are mapped to generators. Clearly
$P_x \fL_G P_x$ is generated by $L_w$ with $ w \in \P(\fgeeplus,
x)$. Let $\T_x$ be the collection of all vectors $\xi_w$, $w=xw$,
such that $w\neq uv$, $u \in \P(\fgeeplus, x)$, $v \in \fgeeplus$.
The subspaces
\[
\V_w = \spn\big\{ A\xi_w : A \in P_x \fL_G P_x \big\} \qfor w \in
\T_x \cup\{\xi_x\}
\]
are mutually orthogonal, reducing for $P_x \fL_G P_x$,  and
satisfy
\[
\sum_{w \in \T_x \cup\{\xi_x\}}\!\!\oplus\,\, \V_w = P_x(\H_G).
\]
For each $w \in  \T_x \cup\{\xi_x \}$ consider the unitary operator
\[
U_w : \V_w \longrightarrow \H_n
\]
defined by $U_w ( \xi_w) = \xi_{\emptyset}$ and
\[
U_w \left (\xi_{u_{i_1}u_{i_2}\dots u_{i_m}w} \right)  = \xi_{i_1 i_2 \dots i_m},
\]
for any $u_{i_1}, u_{i_2}, \dots ,u_{i_m} \in \P(\fgeeplus, x)$.
Then it is evident that a desired unitary operator is given by
$\oplus_{w \in \T_x \cup\{ \xi_x \} } U_w$. \bx

If $ w $ is a loop with source $x \in \V(G)$, then $\| w \| $
denotes the number of  primitive loops whose product equals $w$.

\begin{lem}  \label{restimate}
Let $G$ be a countable directed graph, $x \in \V (G)$ and  let
$e_1,e_2,\dots,e_n \in \E(G)$ with $r(e_i) = s(e_i)=x$,
$i=1,2,\dots,n$, be the distinct loop edges starting at $x$.
Suppose $\rho \in \fM_{G, x}^{w^*}$ satisfies $\sum_{i=1}^{n}\, |
\rho( L_{e_i}) |^2 = r <1$. If $A = P_x A P_x \in \fL_G$ has
Fourier expansion of the form
\[
A \sim \sum_{ \|w \| \geq k} a_w L_w  \qforsome k \geq 1,
\]
then $|\rho(A)| \leq r^k \| A \|$.
\end{lem}

\Prf The proof follows from the identification of $P_x \fL_G P_x$
with $\fL_n$ (Lemma~\ref{secondunitequiv}) and an application of
\cite[Lemma 3.1]{DP1}. \bx

We also need the following lemma, whose proof is similar to that of
Lemma~\ref{secondunitequiv}. (Compare also with \cite[Proposition 3.1]{Sol}).

\begin{lem}  \label{unitequiv}
Let $G$ be a countable directed graph, $x \in \V (G)$ and let $e_1,e_2,\dots,e_n \in \E(G)$
so that $r(e_i) = s(e_i)=x$, $i=1,2,\dots,n$. Let $\alg(P_x , L_{e_1} , L_{e_2}, \dots,L_{e_n})$
be the $w^*$-closed algebra generated by $P_x$ and $L_{e_1} , L_{e_2}, \dots,L_{e_n}$. Then
there exists a $w^*$-bicontinuous
isomorphism from $\alg(P_x ,  L_{e_1} , L_{e_2}, \dots,L_{e_n})$ onto $\fL_n$, which maps
$L_{e_i}$ to $L_i$, $i=1,2,\dots,n$.
\end{lem}

\begin{thm} \label{w*ball}
Let $G$ be a countable directed graph. Let $x \in \V(G)$ and
assume  there exists exactly $n$ distinct loop edges whose source
is $x$. Then $\fM_{G, x}^{w^*}$ equipped with the relative
$w^*$-topology is homeomorphic to $\bbB_{n}^{\circ}$, the open unit
ball of $\bbC^n$.
\end{thm}

\Prf Let $e_1, e_2, \dots ,e_n$, $e_i=xe_ix$, $i= 1,2, \dots,n$, be the
distinct loop edges starting at $x \in \V(G)$ and define
$\phi_x : \fM_{G, x}^{w^*} \longrightarrow \bbC^{n}$
 by the formula
\[
\phi_x(\rho)= (\rho(L_{e_1}), \rho(L_{e_2}), \dots, \rho(L_{e_n})), \quad
\rho \in  \fM_{G, x}^{w^*}.
\]
Since any multiplicative form is completely contractive, the range of
$\phi_x$ is contained in $\bbB_n$. We will show that $\phi_x$ is the desired
homeomorphism.

Arguing as in the proof of Proposition \ref{ball}, we obtaint that
$\phi_x (\fM_{G, x}) \subseteq \bbB_{n}^{\circ}$. Now Lemma~\ref{unitequiv} shows that if
there were $w^*$-continuous multiplicative forms $\rho \in \fM_{G,
x}^{w^*}$ with $\rho \in \vartheta\bbB_n$, then such forms would
also exist on $\fL_n$. But this contradicts~\cite[Theorem
2.3]{DP1} and so $\phi_x (\fM_{G, x}) = \bbB_{n}^{\circ}$.

We now show that $\phi^{-1}_x$ is continuous. Let $A =P_x A P_x
\in \fL_G$ with Fourier expansion $A \sim \sum_{w = xwx} \, a_w
L_w$ and let $\rho, \mu \in \fM_{G, x}^{w^*}$ satisfying $\|
\phi_x(\rho)\|_2 ,\| \phi_x(\mu) \|_2 \leq r <1$. Then given $k
\geq 1$, Lemma~\ref{restimate} shows that
\begin{eqnarray*}
|\rho (A) - \mu(A) |&\leq& \big| (\rho - \mu )\big(\sum_{\|w \|<k}
\, a_w L_w \big) \big| + 2 r^k ||A - \sum_{||w||<k} a_wL_w || \\
&\leq& \big| (\rho - \mu )\big(\sum_{\|w \|<k} \, a_w L_w \big)
\big| + 2 (k+1) r^k ||A ||
\end{eqnarray*}
Notice however that if $w \in \P(\fgeeplus, x)$ contains an edge different from
$e_1,e_2,\dots,e_n$,
then $\rho(L_w)=0$. Hence,
\begin{eqnarray*}
|\rho (A) - \mu(A) |&\leq& \sum_{1\leq i_1 ,, i_2, \dots , i_k
\leq n}
  \, | a_{e_{i_1}e_{i_2}\dots e_{i_k}} (\rho - \mu )(L_{e_{i_1}e_{i_2}\dots e_{i_k}})| \\
  & & + 2(k+1)r^k \| A \|.
\end{eqnarray*}
Thus, since $k$ was arbitrary and $0 \leq r < 1$, the result
follows from a standard approximation argument. \bx

\begin{cor}  \label{w*-vertexcor}
If $G$ is a countable directed graph, then $\fM_{G}^{w^*}$ is a
locally compact Hausdorff space whose connected components
coincide with $\fM_{G, x}^{w^*}$, $x \in V(G)$.
\end{cor}

We now proceed as in Section~\ref{S:quiver}. Given two distinct
vertices $x_1, x_2$ in a directed graph $G$, we define
$\rep_{x_1,x_2}^{w^*}(\fL_G)$ to be the collection of all
$w^*$-continuous, two dimensional nest representations $\pi$ for
which $\rho_{\pi}^{(i)} \in \fM_{G , x_i }^{x*}$, $i=1,2$. Arguing
as in Corollary~\ref{existence}, one shows that
$\rep_{x_1,x_2}^{w^*}(\fL_G) \neq \emptyset$ exactly when there
exists $e \in \E(G)$ such that $e = x_2 e x_1$. Let
\[
\K_{x_1,x_2}^{w^*} =\bigcap\, \{ \ker\pi \mid \pi \in
\rep_{x_1,x_2}^{w^*}(\fL_G) \}
\]
and let $^{\perp}(\K_{x_1,x_2}^{w^*})$ denote the collection of all
$w^*$-continuous functionals on $\fL_G$ that vanish on $\K_{x_1,x_2}^{w^*}$,
equiped with the usual norm. It is a standard result in Functional Analysis
that the dual Banach space of $^{\perp}(\K_{x_1,x_2}^{w^*})$ is isometrically
isomorphic to $\fL_G / \K_{x_1,x_2}^{w^*}$. Therefore
$\fL_G / \K_{x_1,x_2}^{w^*}$ can be equiped with a $w^*$-topology so that
it becomes a dual Banach algebra.

\begin{defn}
Let $\I$ be a $w^*$-closed ideal of a dual
Banach algebra $\B$ and let $\{ B_i\}_{i \in \bbI}$
be a subset of $\I$. The set $\{ B_i\}_{i \in \bbI}$ is said to be a \textit{w*-generating set} for
$\I$ if the $w^*$-closed ideal generated by $\{ B_i\}_{i \in \bbI}$ equals $\I$. (In that case,
the $B_i$'s are said to be the \textit{w*-generators} of $\I$.) The ideal $\I$ is said to be
\textit{$n$-generated} with respect to the $w^*$-topology iff $n$ is the smallest cardinality of a $w^*$-generating
set. If in addition,
the sequential $w^*$-closure, i.e., $w^*$-limits of sequences, of the algebraic ideal generated
by the $n$ generators equals $\I$, then $\I$ is said to be \textit{sequentially $n$-generated}.
\end{defn}

As in Section~\ref{S:quiver}, let $\rad (\fL_G / \K_{x,y}^{w^*})$ be the Jacobson
radical of $\fL_G / \K_{x,y}^{w^*}$. An argument similar to that of Proposition~\ref{radical} shows that
$\rad (\fL_G / \K_{x,y}^{w^*})$ consists of all cosets $A + \K_{x,y}^{w^*} \in \fL_G / \K_{x,y}^{w^*}$
which satisfy $\pi(A)^2 = 0$, for all $\pi \in  rep_{x,y}^{w^*}(\fL_G)$.

\begin{thm}   \label{w*principal}
Let $G$ be a countable directed graph and let $x, y \in \V(G)$ with $x \neq y$.
Assume that there exist
$n$ distinct edges with source $x$ and range $y$. Then the ideal
$\rad (\fL_G / \K_{x,y}^{w^*}) \subseteq \fL_G / \K_{x,y}^{w^*}$ is sequentially n-generated
with respect to the w*-topology of $\fL_G / \K_{x,y}^{w^*}$.
\end{thm}

\Prf Let $e_1, e_2, \dots, e_n$ be the edges with source $x$ and range $y$. An argument similar to
that of Theorem~\ref{principal} shows that the  collection,
\[
L_{e_{i}} + \K_{x,y}^{w^*}, \quad i=1, 2, \dots,n,
\]
is a generating set for $\rad (\fL_G / \K_{x,y}^{w^*})$. Furthermore, the $w^*$-convergence of the Cesaro-type sums
shows that $\rad (\fL_G / \K_{x,y}^{w^*})$ is sequentially n-generated.

By way of contradiction assume that there exists a generating set
\[
A_j + \K_{x,y}^{w^*}, \quad j=1, 2, \dots,m,
\]
for $\rad (\fL_G / \K_{x,y}^{w^*})$ with $m<n$. By Lemma~\ref{reduction}, we may assume that
$A_j = P_y A_j P_x $. Therefore, there exist scalars $a_{j, i}$ so that
\[
\Phi_1(A_j)= a_{j,1}L_{e_{1}} +a_{j,2}L_{e_{2}}+\dots+a_{j,n}L_{e_{n}} ,
\]
where $\Phi_1$ is the contractive idempotent on $\fL_G$ defined by the formula
\[
\Phi_1\Big(\sum_w \, a_wL_w\Big)=\sum_{|w|=1} \, a_wL_w.
\]
Let $\M$ be the subspace of $\fL_G$ generated by $L_{e_{i}}$, $i=1,2,\dots,n$, and let
$\M_0 \subset \M$ be the linear span of
\[
a_{j,1}L_{e_{1}} +a_{j,2}L_{e_{2}}+\dots+a_{j,n}L_{e_{n}}, \quad j=1,2,\dots,m.
\]
Let $X \in \M \backslash \M_0$. Clearly, $X+ \K_{x,y}^{w^*}$ belongs to $\rad (\fL_G / \K_{x,y}^{w^*})$.
Since the collection $A_j + \K_{x,y}^{w^*}$, $j=1, 2, \dots,m$, is generating
for $\rad (\fL_G / \K_{x,y}^{w^*})$, the coset $X+ \K_{x,y}^{w^*}$ can be approximated in the $w^*$-topology
by a sequence $\{ X_n + \K_{x,y}^{w^*} \}_{n=1}^{\infty}$ consisting of finite sums of the form
\[
X_n + \K_{x,y}^{w^*} = \sum_{k}\, B_{k}^{(n)} \hat{A}_{k}^{(n)} C_{k}^{(n)} + \K_{x,y}^{w^*},
\]
where, $B_{k}^{(n)} , C_{k}^{(n)} \in \fL_G$, and $\hat{A}_{k}^{(n)} \in
\{A_1, A_2, \dots A_m \}$. Since $\{ X_n + \K_{x,y}^{w^*} \}_{n=1}^{\infty}$ is $w^*$-convergent,
it is bounded. Hence we can choose $D_n \in \K_{x,y}^{w^*}$, $n \in \bbN$, so that
the sequence $\{ X_n + D_n \}_{n=1}^{\infty}$ is bounded as well. Passing
to a subsequence if necessary, we may assume that the sequence $\{ X_n + D_n \}_{n=1}^{\infty}$
is $w^*$-convergent to $X+D$, for some $D \in \K_{x,y}^{w^*}$.
Lemma \ref{Phi1}
shows now that
$\Phi_1(P_y (D_n - D) P_x)= 0$, for all $n \in \bbN$. It is also easy to see that
\[
\Phi_1 (B A_j C) = \mu_j \big( a_{j,1}L_{e_{1}} +a_{j,2}L_{e_{2}}
+\dots+a_{j,n}L_{e_{n}}\big), \quad \mu_j \in \bbC,
\]
for all $j= 1,2,\dots,m$ and $B,C \in \fL_G$. Hence,
\begin{equation*}
\Phi_1 \left(  P_y \left( \sum_{k}\, B_{k}^{(n)} \hat{A}_{k}^{(n)}
C_{k}^{(n)}\, +\, D_n -D \right)P_x \right) \in \M_0.
\end{equation*}
A simple approximation argument with the contractive map $\Phi_1$ shows now
that $\Phi_1(X) = X$ belongs to $ \M_0$, which is a contradiction.
\bx

A verbatim
repetition of the proof of Theorem~\ref{main quiver} proves now the
following.

\begin{thm} \label{main freesem}
Let $G$, $G'$ be countable directed graphs. Then there exists a
$w^*$-bicontinuous isomorphism $\tau : \fL_G \longrightarrow
\fL_{G'}$ if and only if the graphs $G$ and $G'$ are isomorphic.
\end{thm}

\begin{defn}
If $\A$ (resp. $\B$) is an algebra on a Hilbert space $\H_{\A}$ (resp. $\H_{\B}$)
then $\A$ and $\B$ are said to be \textit{spatially isomorphic} if there exists an
invertible operator $S : \H_{\A} \longrightarrow \H_{\B}$ so that $S \A S^{-1} = \B$.
If $S$ is a unitary operator then $\A$ and $\B$ are said to be \textit{unitarily isomorphic}.
\end{defn}

\begin{cor}
Two free semigroupoid algebras $\fL_G$ and $\fL_{G'}$ are spatially isomorphic
if and only if they are unitarily isomorphic.
\end{cor}

\Prf A similarity between two free semigroupoid algebras
$\fL_G$ and $\fL_{G'}$ is a map
satisfying the requirements of Theorem~\ref{main freesem}. Therefore, the graphs $G$ and $G'$
are isomorphic and so $\fL_G$ and $\fL_{G'}$ are unitarily isomorphic \cite[Theorem 9.1]{KrP}.
The other direction is trivial.
\bx

\begin{rem}
It should be emphasized that a spatial isomorphism between two free
semigroupoid algebras $\fL_G$ and $\fL_{G'}$ need not transform $\lambda_G$ into
$\lambda_{G'}$. The relation between the spatial isomorphism and the unitary isomorphism guaranteed
by our analysis can be quite mysterious. This is partly due to the fact that our understanding
of the structure of the automorphism group of a free semigroup algebra is still limited.
\end{rem}

In order to prove that the graph forms an invariant for algebraic isomorphisms of
free semigroupoid algebras, it remains to show that algebraic isomorphisms between
such algebras are $w^*$-continuous. In what follows we do this for a special class of graphs. The
general case remains open.

\begin{thm}  \label{automatic cntns LG}
Let  $G$ be a countable directed graph with no sinks and let $\A$
be any Banach algebra. If $\tau : \A \longrightarrow \fL_{G}$ is
an epimorphism of algebras then $\tau$ is automatically
continuous.
\end{thm}

\Prf The proof is identical to that of Theorem \ref{automatic cntns}.
\bx

One of the main results in \cite{JK1} asserts that the algebra $\fL_{G}$
satisfies the property $\bbA_1$, provided that $G$ has no sinks. Therefore for the
free semigroupoid algebras considered in the rest of this section, the WOT and $w^*$-topology
coincide.

\begin{lem} \label{G0k}
Let $G$, $G^{\prime}$ be locally finite directed graphs with no
sinks. If $\tau : \fL_G  \longrightarrow \fL_{G'}$ is an algebraic
isomorphism, then $\tau (\fL_G^{0,k})$, for $k\geq 1$, is
$w^*$-closed.
\end{lem}

\Prf Let $\fL_{G, x}^{0,k}$ be the collection of all $A \in \fL_G$ which can be written
as a sum
\[
A = (I-P_x)A_1  +  P_x A_2,
\]
where $A_2 \sim \sum_{|w| \geq k}a_w L_w$. It is easy to verify that $\fL_{G, x}^{0,k}$ is
$\wot$-closed and that $\fL_{G}^{0,k} = \bigcap_{x \in \V(G)} \fL_{G, x}^{0,k}$. Therefore it
suffices to show that for each $x \in \V(G)$, the set $\tau(\fL_{G, x}^{0,k})$ is $\wot$-closed.

Let $\{ B_{\alpha}\}_{\alpha \in \bbA}$ be a net in $\tau (\fL_{G, x}^{0,k})$ converging
to $B$.  The Krein-Smulian Theorem implies that there is no loss of generality
assuming that $\{ B_{\alpha}\}_{\alpha \in \bbA}$ is bounded in norm.
Let $\{ A_{\alpha} \}_{\alpha \in \bbA}$ be a bounded net in $\fL_G$ so that $B_\alpha = \tau
(A_\alpha)$, for all $\alpha \in \bbA$ (recall that $\tau^{-1}$ is continuous).
Corollary 4.8 in~\cite{JK1} implies that $A_\alpha$ can be written as a sum
\begin{equation} \label{alpha}
A_\alpha = (I-P_x)A^{(\alpha)}_{1} \, + \, \sum_{|w| = k }L_{xw}A_{w}^{(\alpha)} ,
\end{equation}
where $ A^{(\alpha)}_{1} , A_{w}^{(\alpha)} \in \fL_G$, for all $\alpha \in \bbA$
and $w \in \fgeeplus$. (Note that the local finiteness of $G$ guarantees that
the sum in (\ref{alpha}) is finite.)
Since $\{ A_{\alpha} \}_{\alpha \in \bbA}$ is a bounded net, the nets $\{ A^{(\alpha)}_{1}
\}_{\alpha \in \bbA}$
and $\{ A_{w}^{(\alpha)} \}_{\alpha \in \bbA}$ are also bounded. Since $\tau$ is continuous,
the nets $\{\tau( A^{(\alpha)}_{1} ) \}_{\alpha \in \bbA}$
and $\{ \tau( A_{w}^{(\alpha)}) \}_{\alpha \in \bbA}$ are bounded as well.
The compactness of the closed ball in the $w^*$-topology, combined with a diagonal argument,
shows that the net
$\{ B_{\alpha}\}_{\alpha \in \bbA}$ has a subnet converging to an element of the form
\[
 \tau (I-P_x) \tau (A_{1}) \, + \, \sum_{|w| = k }\tau (L_{xw}) \tau (A_{w})
\]
and so $B = \tau \left((I-P_x)A_{1} \, + \, \sum_{|w| = k }L_{xw}A_{w} \right) \in  \fL_{G,
x}^{0,k}$,
as desired.
\bx

\begin{lem}  \label{strong*}
Let $G, G'$ be countable directed graphs with no sinks and let
$\tau : \fL_G \longrightarrow \fL_{G'}$ be an algebraic
isomorphism. Then $\sum_{x \in \V(G) }\, \tau (P_x)= I$ in the
$\mbox{strong}^{\,*}$ topology.
\end{lem}

\Prf Let $\A$ be the unital $\ca$-algebra generated by $P_x$, $x \in \V(G)$.
Since $\A$ is abelian and $\tau$ can be thought as a representation of $\A$, there exists an
invertible
operator $S \in B ( \H_{G^{'}})$ so that the mapping $\hat{\tau} : \fL_G \longrightarrow \B(
\H_{G^{'}})$,
defined as $\hat{ \tau} (A)= S \tau (A) S^{-1}$, $ A \in \fL_G$, is a $*$-representation on $\A$.
It suffices to show that $\sum_{x \in \V(G) }\, \hat{ \tau }(P_x)= I$, in the strong topology.

Note that since the projections $P_x$, $x \in \V(G)$, are
mutually orthogonal, the same is true for $\hat{ \tau}(P_x)$, $x \in \V(G)$.
Consider
the projection
\[
Q \, = \,\sum_{x \in \V(G) } \, \hat{ \tau}(P_x) \in \hat{ \tau}(\fL_{G}).
\]
Then $(I-Q)\hat{ \tau}(P_{x}) = 0$, and so
$\hat{ \tau}^{-1}(I-Q)P_{x}=0$, for all $x \in \V(G)$.
Since, $\sum_{x \in \V(G) } \,P_{x} = I$, we obtain that $\hat{ \tau}^{-1}(I-Q)=0$.
Therefore, $\sum_{x \in \V(G) }\, \tau (P_x)= I$, as desired.
\bx

The proof of the following theorem is modelled on the proof from
\cite{DP1} for the special case of $\fL_n$.

\begin{thm}\label{wotcont}
Let $G, G'$ be locally finite directed graphs with no
sinks. Then every algebraic isomorphism $\tau : \fL_G
\longrightarrow \fL_{G'}$ is $w^*$-continuous.
\end{thm}

\Prf In \cite{JK1} it was shown that the $w^*$ and $\wot$ topologies
on $\fL_G$ coincide when $G$ has no sinks. Thus, by an application
of the Krein-Smulian Theorem, it follows that $\tau$ is
$w^*$-continuous if and only if $\tau$ is $\wot$-continuous on
every closed ball of $\fL_G$.

Let $A_\alpha$ be a bounded net of operators in $\fL_G$ which
converge $\wot$ to zero.
Let $x_1, x_2, \ldots$ be an enumeration of $G$ and let
$\xi_\phi = \sum_{n \in \bbN} \, 2^{-n}\xi_{x_n}$.
By an elementary argument (see for instance \cite[Lemma 4.4]{DP1}), it suffices to show that
\[
\lim_\alpha \langle \tau(A_\alpha) \xi_{\phi}, \zeta \rangle =0,
\]
for all $\zeta$ in a dense subset of $\H_{G'}$.
We shall obtain a distinguished dense subset as follows. From Lemma \ref{G0k} we know that
\[
\J_k \equiv \tau (\fL_G^{0,k})
\]
is $\wot$-closed.
Observe that $\cap_{k\geq 1} \J_k = \{ 0\}$ since
\[
\tau^{-1} (\cap_{k\geq 1} \J_k) \subseteq \tau^{-1}(\J_k) =
\fL_G^{0,k} \qfor k\geq 1.
\]
It follows from \cite[Theorem~5.2]{JK1} that
\[
\cap_{k\geq 1} \ol{\J_k\H_{G'}} = \{ 0\},
\]
and hence $\cup_{k\geq 1} (\J_k\H_{G'})^\perp$ is dense inside
$\H_{G'}$.

Now let $\zeta\in (\J_k \H_{G'})^\perp$ and let $\epsilon > 0$.
Using Lemma~\ref{strong*} we can choose a subset $\S
\subseteq \V(G)$ so that $\V(G) \backslash \S$ is finite and
\[
\| \sum_{x \in \S}\,\tau( P_x )^{*}\zeta
\| \leq \epsilon.
\]
Therefore, it suffices to identify $\alpha_{0}$
such that
\begin{equation} \label{inS}
\left| \langle \tau (P_x A_\alpha ) \xi_{\phi} , \zeta \rangle \right| \leq \frac{\epsilon}{| \S |},
\end{equation}
for all $\alpha \geq \alpha_0$ and $ x \in \V(G) \backslash \S$.
Suppose
$ A_\alpha \sim \sum_w a_{w}^{(\alpha)} L_{w}$ and decompose $P_x A_\alpha =
B_{\alpha ,x} + C_{\alpha ,x}$ where
\[
B_{\alpha ,x} = \sum_{|w|<k} a_{xw}^{(\alpha)} L_{xw}
\]
and
\[
C_{\alpha ,x}= A_{\alpha ,x} - B_{\alpha ,x} = \sum_{|w| = k} L_w A_w^{(\alpha)}
 \in \fL_{G'}^{0,k}.
\]
By construction, $\langle \tau(C_\alpha) \xi_\phi, \zeta \rangle =0$  since $\tau(C_\alpha)
\xi_\phi$ belongs to $\J_k
\H_{G'}$. On the other hand, $a_w^{(\alpha)} = \langle A_\alpha \xi_{s(w)}, \xi_w \rangle$, so that
$\lim_\alpha a_w^{(\alpha)} =0$ for all $w\in \fgeeplus$, $|w| < k$. Since
the sum $ B_{\alpha ,x} = \sum_{|w|<k} a_{xw}^{(\alpha)} L_{xw}$ is finite ($G$ is locally finite), the desired
inequality (\ref{inS}) follows.
\bx

\begin{cor}
Let $G$, $G'$ be locally finite directed graphs
with no sinks. Then $\fL_G$ and
$ \fL_{G'}$ are
isomorphic as algebras if and only if the graphs $G$ and $G'$ are isomorphic.
\end{cor}

\vspace{0.1in}

{\noindent}{\it Acknowledgements.} We would like to thank the
referee for several constructive suggestions on the initial draft and for
bringing to our attention the work in \cite{FR, FMR}.
The first author was partially
supported by a research grant from ECU and the second author
by an NSERC research grant and start up funds
from the University of Guelph. We thank David
Pitts for enlightening conversations and Alex Kumjian for helpful
comments on the literature.

\end{document}